# Research on High-precision Detection Technology for Underground Space Information Data

**Zuo Yanhong, Zhou Chao, Xia Shilong, Yang kun**

*School of Mechanical and Electrical Engineering, Anhui Jianzhu University, Hefe 230601, China;*

*\* Correspondence: zuoyh626@sohu.com*

**Abstract**: The quality of underground space information data has become a major problem endangering the safety of underground space. After research and analysis, we found that the current high accuracy information data remote detection methods are limited to the detection of overground spaces objects, and are not applicable to the detection of various information data in underground space. In this paper, we analyze the spectral properties of the fractional-order differential (FDO) operator, and establish mathematical model of remote transmission and high-precision detection of information data, which realizes the functions of high-precision and remote detection of information data. By fusing the information data to detect the mathematical model in a long distance and with high accuracy, A mathematical model has been established to improve the quality of underground spatial information data. Through the application in engineering practice, the effectiveness of this method in underground space information data detection is verified.

**Keywords:** Underground space. Information detection. Fractional differentiation. High accuracy data

# 1 INTRODUCTION

Currently, the development of larger and more complex underground spaces has become a trend for future underground space development [1]. Underground space has the characteristics of "deep, large, clustered, and hidden". Once an accident occurs, it can easily cause a catastrophic impact [2-3]. However, the existing safety perception and prediction technology for underground space structures cannot guarantee the safety of underground space structure construction and operation. In particular, urgent problems still need to be solved in terms of state perception, disease identification, and mechanistic analysis [4]. In addition, the working environment of underground space measuring instrumentFOD is complex, spatio-temporal monitoring data are unusually diverse, and the quality is uneven, which often causes problems such as "unclear detection, inaccurate measurement, and inaccurate judgment" of the underground space, thus leading to an incorrect estimation of the safety status of the underground space, which may endanger the safety of the underground space [5]. Therefore, improving the global awareness of underground spaces and the quality of data detection has become an urgent problem to be solved in the construction of underground spaces.

The detection objects of underground space information include the detection of underground tangible objects and intangible material parameters; however, the essence of realizing their information data detection quality lies in improving the detection accuracy of information data and the strength of signal transmission. The tangible objects in the underground space include the underground space structure, object contour shape, production equipment position, and moving object trajectory. At present, the detection methods of underground space are mainly based on ground-penetrating radar [6], and the portability and high-resolution characteristics of GPR [7-8] make GPR widely used in subgrade surveys and underground space pipeline detection [9]. However, the penetration depth of ground-penetrating radar is limited, and the effective detection depth generally does not exceed 10m [10-11], it can only detect shallow urban facilities and has no detection capability for underground targets below 30m. The existing shallow surface wave method [12-14] has difficulty meeting the detection requirements of urban geophysical exploration owing to its weak anti-jamming capability. Reference [15-16] proposed a load differential radiation pulse on a transient electromagnetic high-property radiation source for pulse scanning detection to solve the problems of urban electromagnetic interference and insufficient harmonic components emitted by radiation sources. However, the detection method of pulse scanning alone is incomplete and a corresponding imaging method is required. More importantly, the above methods can only be applied to the detection of information tens of meters below the ground, and they cannot be applied to underground spaces hundreds of meters deep.

The objects of detection of intangible substances in underground spaces include gas, pressure, temperature, humidity, noise, and stress. At present, the detection methods for the parameters of underground intangible substances are still the same as those for intangible substances in above-ground space. Reference [17-19] designed a hardware system based on radar and realized the real-time detection function of underground space related information by enlarging the detection information. However, the measurement error cannot be eliminated and applied to underground spaces with tortuous paths. Reference [20] proposed a three-frequency resonance transmission scheme to solve the problems of low efficiency in the conventional single-frequency transmission method and high voltage stress in the multifrequency pseudorandom transmission method. However, the scheme only solves the problem of signal transmission efficiency and does not improve the quality of information data. Reference [21-22]The self-developed rotary drilling system and three-dimensional flexible boundary loading device were developed to carry out the rotary drilling model and realize the detection technology of stress betteen rocks in underground space. However, they have limitations in the application space and cannot be used for remote detection of information data. ZhengXuezhao et al.[23] studied a multi-information drilling detection device based on key technologies, such as multimedia information collection, synchronous transmission, and underground explosion-proof facilities, and realized the collection function of key information data of underground space. However, this method is based on drilling from the ground to the underground space and placing the equipment in the underground space for later detection, so there are shortcomings of low efficiency and information lag. Bin Sun et al. [24] proposed a bionic artificial intelligence algorithm driven by temperature data to detect a fire source in the three-dimensional space of an underground pipe gallery. However, this method can be used only for fire prevention in underground pipelines.
More importantly, the above methods have the following shortcomings in the information collection of underground spaces.
  (1) They are mainly used to detect the position and shape of tangible objects, and there are shortcomings in the field of high-precision remote detection that cannot be applied to object property parameters. Therefore, there are limitations in the application field.
  (2) Their essence is to highlight the difference between information data to improve the detection accuracy of the difference data, but they fails to eliminate the detection error caused by various factors in the process of long-distance transmission of in-formation data fundamentally.
  (3) In the detection process, because the detection object is in a random transformation state, the above information detection system does not have the ability to adjust the technical parameters with the changes in the tested object in a timely manner. Therefore, it affects the accuracy and real-time detection of the information.

(4) In the detection of underground space information data, the above methods have the disadvantage of large energy loss in application, which restricts the depth and data accuracy of detection information; therefore, there are limitations in the application space.

Therefore, to date, we have not yet found any detection methods that can effectively improve the quality of information data for various types of detection objects in underground spaces.. To solve these problems, our team has engaged in research on methods using fractional-order calculus theory in data processing for many years and found that the FOD operator has the dual function of improving the signal strength and reducing the variability between information data. Therefore, it can effectively improve the transmission distance and detection accuracy of information [25-32]. Based on previous research, in this study, we applied fractional calculus theory to the detection of underground space information data.

## 2 FRACTIONAL CALCULUS THEORY

At present, the commonly used definitions are those of Grunwald–Letnikov (G-L), Caputo, and Riemann–Liouville (R-L) [33]. Because the G-L definition has the advantage of fast calculation speeds, it is widely used in the engineering field [34]. Therefore, the G-L definition was used in this work to study the application of fractional calculus theory for detection data processing.

### 2.1 Fractional-order Calculus Definition

The definition of G-L is derived from the extension of the order of the classical integer order differential of continuous function from the integer order to the fraction, and the calculation of the limit of the original integer order differential measurement difference myopia recurrence formula. If $R$ is used to represent the real number field, and $[v]$ is used to represent the integer part of v, assuming that the signal $f(t)$ is satisfied $s(t) \in [a, t]$ ($a < t, a \in R, t \in R$). If $Z$ is used to represent the integer field, let $P \in Z$, $s(t)$ be used to make the signal continuous derivative at $P+1$ order. When $v > 0$ is used, $P$ at least takes $[v]$, then the $v$-order derivative of the signal is defined as follows.

$$_a D_t^v f(t) = \lim_{h \to 0} f_h^{(v)}(t) = \lim_{h \to 0} h^{-v} \sum_{j=0}^{\left[\frac{b-a}{h}\right]} (-1)^j \binom{v}{j} f(t-jh). \tag{1}$$

In (1), $0 \leq n-1 < v < n$, and $\binom{v}{j}$ is a binomial coefficient, defined as

$$\binom{v}{j} = \frac{v(v-1)(v-2)....(v-j+1)}{j!}. \tag{2}$$

where $h = (b-a)/n$, and $n = [(b-a)/h]$. When $h \to 0$, $n \to \infty$.

### 2.2 Spectral Characteristics of Fractional Calculus Operators

For an arbitrary square productable energy signal $f(t) \in L^2(R)$, it is assumed that its Fourier transform is $F(\omega)$, then the Fourier transform of the $v$-order FOD $D^v f(t)$ can be obtained as follows:

$$\frac{d^v f(t)}{dt^v} \Leftrightarrow (i\omega)^v F(\omega) \tag{3}$$

The characteristic function of the signal $v$-order FOD operator $D^v$ is

$$D(\omega) = \alpha(\omega)\theta(\omega) = (i\omega)^v \tag{4}$$

where $\alpha(\omega) = |\omega|^v$ and $\theta(\omega) = \frac{\pi v}{2}\text{sgn}(\omega)$

Based on the filter function, we can draw the spectral characteristic curves of the FOD operator in signal processing, as shown in Figures 1 respectively. After the analysis of the spectral characteristic curves, we can obtain the following characteristics of fractional-order calculus for signal processing.

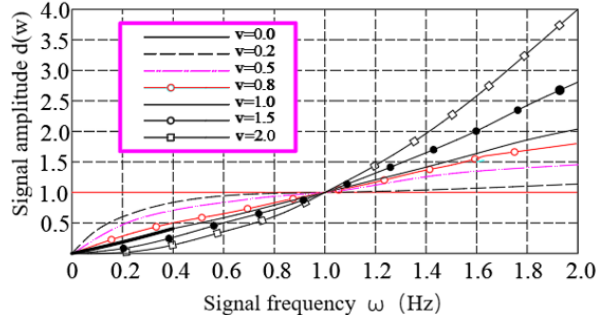

Fig. 1. Amplitude–frequency characteristic curves of FOD operator

(1) It has the function of non-linearly retaining the very -low -frequency components of the various signals while boosting the high-frequency components of the signal. This can effectively enhance the middle- and high-frequency parts of the signal, and the amplitude tends to increase non-linearly and rapidly with the increase in the frequency and fractional order of the derivative.

(2) When ω > 1, with the increase of fractional order v and signal frequency ω, the variability between the signal enhancement coefficients of the FOD operators at different orders tends to decrease, and with the increase in the frequency, the enhancement effect of the differential operators of the same order on the signal intensity at different frequencies is basically the same.

## 2.3 Application of Fractional Calculus

In the expression of the fractional-order calculus under the G-L definition, the time $t$ in the signal $f(t)$ is the impact parameter of the signal amplitude, where $t \in [t_1, t_2]$. If $t_i = x_i$, and $[t_1, t_2] \Leftrightarrow [a, b]$, reference (1), we can obtain the equation as follows:

$$_a^G D_t^v f(t) \Leftrightarrow \, _a^G D_x^v f(x) = \lim_{h \to 0} h^{-v} \sum_{i=0}^{\left[\frac{b-a}{h}\right]} (-1)^i \binom{v}{i} f(x - ih). \tag{5}$$

In (5), $0 < n - 1 < n$, and

$$\binom{v}{i} = \frac{v(v-1)(v-2)\ldots(v-i+1)}{i!}. \tag{6}$$

In (5), where $h \to 0$, $n \to \infty$, and $h = (b-a)/n$.

In recent years, many experts and scholars have studied the theory of fractional-order calculus and found that it is suitable for the study of signals with undesirable characteristics, such as nonlinearity, noncausality, and non-stationarity, and for various data applications [31]. Various research results have been used for technical problems in the study of temperature field distributions, image processing, mechanical analysis, and detection technology [35-39]. For example, WANG Bao et al. [40] applied fractional partial derivatives for the thickness design of high-temperature protective clothing under actual limited conditions . Zhou et al. [41] applied fractional derivatives to prove their advantages in image denoising and reducing the step effect, as well as in signal denoising and super-resolution reconstruction. SHEN Tianlong [42] applied fractional partial derivatives in the field of fluid mechanics. The working environment of underground space is complex and changeable, all kinds of information have the characteristics of irregular change. Therefore, FOD operators can be used in the processing of various underground signals. It can solves the problem of application field limitation in current underground space information detection methods successfully. High-precision detection

## 3 HIGH - PRECISION DETECTION METHOD OF UNDER-GROUND SPACE INFORMATION BASED ON FRACTIONAL DIFFERENTIAL ALGORITHM

### 3.1 Fundamental Principle

According to Figure 1 and Formula (1), the essence of fractional calculus processing of signal $F(t)$ is to draw its spectrum characteristic curve by means of broken line fitting, where $t \in (t_1, t_2)$. In the process of information detection, information data is affected by various factors such as working environment and equipment performance, resulting in unpredictable detection errors. In equation (5), if we take x as the parameter value of a certain influencing factor, we can obtain the functional formula $F(x)$ between the signal measurement value and the influencing factor, where $x \in (a, b)$. If h is the step value of $x$ changing between $(a, b)$, the degree of fitting required is $n = [b − a]/h$. According to (2) and Figure 1, when the value of $h$ is different, the spectral curve of fractional calculus defined by G-L is shown in Figure 2.

In Figure 2, if the step values are set to $h$, $2h$, $3h$, and $4h$ respectively, draw information fitting polylines based on the definition of fractional order differential G-L. By comparing the fitting line with the real line, it can be seen that the detection error $S$ is proportional to the step value $h$. Therefore, the information data detection accuracy $S$ required by the detection system can be obtained by adjusting the step value $h$.

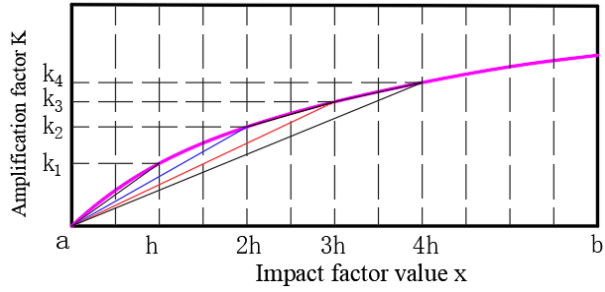

Fig. 2. Data fitting accuracy for FOD operators of the same order

### 3.2 Mathematical Model

According to the signal characteristics and the spectral properties of the FOD operator, the optimal fractional order $v$ for the FOD processing of the data is selected. Referring to equation (5), the system obtains the functional equation $F(x)$ that relates the data value $F(x_i)$ and the impact parameter $x_i$, where $x \in (a, b)$. According to (4), the derivative processing equation of $F(x)$ at fractional order $v$ can be obtained as follows:

$$_aD_x^v F(x) = \lim_{h \to 0} F_h^v(x) = h^{-v} \sum_{i=0}^{n} (-1)^i \binom{v}{i} F(x_i - ih), \qquad (7)$$

where $\binom{v}{i} = \dfrac{v(v-1)(v-2)\ldots(v-i+1)}{i!}$, and $n = [b-a]/h$

According to the above analysis of the spectral characteristics of the FOD operator, the FOD operator has the effect of enhancing the signal intensity. Suppose that when $x$ takes the value interval $[a, b]$, the amplification factor of the data value of the function $F(x)$ after differentiation of order $v$ is $K$. For comparison purposes, the data value after $v$-order differentiation is divided by the amplification coefficient $K$. The result is combined with the formula for the standard deviation $S$ of the data, as follows:

$$S = \left( \frac{(F_a^v - \bar{F}^v)^2 + (F_{a+h}^v - \bar{F}^v)^2 + (F_{a+2h}^v - \bar{F}^v)^2 \ldots + F(F_b^v - \bar{F}^v)^2}{nK} \right)^{0.5} = \left( \frac{\sum_{i=0}^{n_1}(F_{a+ih}^v - \bar{F}^v)^2}{nK} \right)^{0.5}, \qquad (8)$$

where

$$F_{a+ih}^v = \lim_{h \to 0} F_h^v(a+ih) = h^{-v} \sum_{i=0}^{n}(-1)^i \begin{bmatrix} v \\ i \end{bmatrix} F(x_i), \qquad (9)$$

$$\bar{F}^v = \frac{F_a^v + F_{a+h}^v + F_{a+2h}^v + \text{L} + F_b^v}{n} = \sum_{i=0}^{n} F_{a+ih}^v / n, \qquad (10)$$

$$K = \frac{F_a^v + F_{a+h}^v + F_{a+2h}^v + \text{L} + F_b^v}{F_a + F_{a+h} + F_{a+2h} + \text{L} + F_b} = \sum_{i=0}^{n} F_{a+ih}^v \Big/ \sum_{i=0}^{n} F_{a+ih} . \qquad (11)$$

By substituting (8) and (9) into (7), we obtain

$$S = \left( \frac{\sum_{i=0}^{n}(F_{a+ih}^v - \bar{F}^v)^2}{n \sum_{i=0}^{n} F_{a+ih}^v \Big/ \sum_{i=0}^{n} F_{a+ih}} \right)^{0.5}. \qquad (12)$$

If the accuracy of the data is set based on $S_g$, then $S$ should satisfy

$$\begin{cases} S \leqslant S_g \\ S_g - S \approx 0 \end{cases}. \tag{13}$$

Initially, the step value $h$ is chosen. Combined with the known parameter values $a$ and $b$, we can calculate the inter-data accuracy for the step value $h$.

$$S = \left( \sum_{i=0}^{[b-a]/h} (F_{a+ih}^v - \overline{F}^v) / nK \right)^{0.5}. \tag{14}$$

The error between the standard deviation $S$ and $S_g$ is compared, and the step value $h$ is gradually adjusted until the conditions shown in (12) are met. The corresponding step value $h$ is the best step value to meet the set accuracy of the detection system.

### 3.3 Implementation Steps

Based on the above undamental principle and mathematical model, the implementation steps of high-precision data detection method for underground space are as follows:
(1) Based on the characteristics of the detected signal and the spectral characteristics of the FOD operator, the optimal fractional derivative order v is selected for the data F(x).
(2) Referring to (1) and (7), the derivative-processed model Fv(x) of the corresponding data function F(x) is established at the fractional order v.
(3) The functional equation S(h) relating the standard deviation S of the data and the step size h is established based on the calculated value of Fv (x).
(4) The initial step value h is set, and use MATLAB mathematical simulation software to calculate the standard deviation S (h) of the standard deviation S(h) of the detection data when the step value is h is calculated based on the functional equation S(h).
(5) The error between the standard deviation S(h) and the specified standard deviation Sg threshold of the underground space information detection system are compared, and the step value h is continuously adjusted based on the error value between them.
(6) If the standard deviation S(hk) is slightly less than the specified standard deviation Sg threshold, hk is selected as the best step value to meet the accuracy requirement of the underground space information detection system.

If the standard deviation $S(h)$ is slightly less than the specified standard deviation $S_g$ threshold, $h$ is selected as the best step value to meet the accuracy requirement of the information detection system.

## 4 LONG-DISTANCE TRANSMISSION METHOD OF UNDER-GROUND SPACE INFORMATION BASED ON FRAC-TIONAL DIFFERENTIAL ALGORITHM

### 4.1 Fundamental Principle

By combining the amplitude–frequency characteristics of the FOD operator shown in Figure 1 and the FOD processing equation under the G-L definition shown in (1), suppose the system detects an energy signal $G(t)$, where $t \in (t_1, t_2)$. Referring to (5), we can get the corresponding equation $G(x)$, $x \in (a, b)$. Its local enlargement of the spectral characteristics of the function $G(x)$ in the interval $x \in (a, b)$ after differential processing of different orders is shown in Figure 3. When the fractional order $v_i$ and the impact parameter $x_i$ are the same, the amplification coefficient $K_i$ of the FOD operator is a fixed value. If the energy data $G(t)$ is processed by the differential operator of fractional order $v_i$, and the amplification coefficient of the data is $k_i$, then the amplification coefficient of the signal after $m$ cycles of applying the differential operator of fractional order $v$ is denoted as $K = k_i^m$, and the data after processing is denoted as $G_n(x) = k_i^m \overline{G}_h^v(x)$. Therefore, the detection system can realize the long-distance detection function of data by adjusting the number of cycles of the FOD processing of the detection data to $m$.

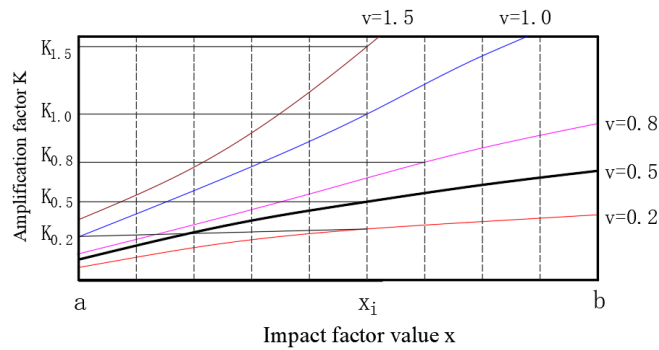

Fig. 3. Differential operator signal characteristics under different orders

## 4.2 Mathematical Model

Since the interval of values of $x$ in the function $G(x)$ is $[a, b]$, where $h$ is the step value, i.e., $x_i = x_{i-1} + h$, and the corresponding step number $n = [(b - a)/h]$. Thus, we can obtain the following equation.

$$\bar{G}(x) = \sum_{i=0}^{n} G(a+ih)/n \cdot \tag{15}$$

The functional equation of the function $G(x)$ after differentiation of order $v$ is $G^v(x)$. From the FOD in (1) and (2), the derivative of order $v$ for $G(x)$ can be obtained as follows:

$$_aD_x^v G(x) = \lim_{h \to 0} G_h^{(v)}(x) = h^{-v} \sum_{i=0}^{n} (-1)^i \binom{v}{i} G(x-ih), \tag{16}$$

where $\binom{v}{i} = \dfrac{v(v-1)(v-2)\ldots(v-i+1)}{i!}$.

Therefore, the average value of the function $G(x)$ processed by the differentiation operator of order $v$ is

$$\bar{G}^v(x) = \sum_{i=1}^{n} G^v(x_i)/n = \frac{h^{-v}}{n} \sum_{i=0}^{n}\sum_{i=0}^{n} (-1)^i \binom{v}{i} G(x-ih), \tag{17}$$

According to the formula of the amplification coefficient $K$, the amplification coefficient $k$ of the function $G(x)$ in the interval $[a, b]$ after the function value is processed by the differentiation operator of order $v$ can be obtained as

$$K = \bar{G}^v(x)/\bar{G}(x) = \frac{h^{-v}}{n} \sum_{i=0}^{n}\sum_{i=0}^{n} (-1)^i \binom{v}{i} G(x_i) / \sum_{i=0}^{n} G(a+ih) \cdot \tag{18}$$

According to the working principle of the long-distance transmission of the data above, it is known that the signal strength will be amplified $k$ times each time the data function $G(t)$ is processed by the $v$-order differential operator. Therefore, if the detection system needs to achieve the required amplification factor $K_g$ for the set long-distance transmission targets under the existing conditions, it must satisfy $K_g \leq k^m$, where $m$ is the number of FOD required to achieve the set detection target, so the required number of $v$-order differential operator processes $m$ should satisfy

$$K^m = \left(\bar{G}^v(x)/\bar{G}(x)\right)^m \geq K_g \cdot \tag{19}$$

The number of iterations directly affects the efficiency of the system, so in the process of increasing the data intensity, the number of fractional-order differentiations is set to

$$K_g \leq K^m < K^{m+1}. \tag{20}$$

where $m$ is the number of cycles in the differential operator for the data.

## 4.3 Implementation Steps

The premise of realizing remote detection of underground space information data is to obtain the detection value $G(x_i)$ of detection data $G(t)$ at different time points in real time, and then analyze the main impact parameter $x$ and its range of values $[a, b]$ based on the data. This lays the foundation for the long-distance transmission of the detection data. The specific implementation process of the long-distance detection method of underground space information data is as follows.

(1) The detection values G(xi) of the signal G(t) of underground space information at different time points are acquired, and the impact parameter x and its value interval [a. b] are analyzed.
(2) The acquired detection values G(xi) and the values of their corresponding impact parameter xi are used to fit a function G(x).
(3) Considering the signal characteristics of G(t) and the spectral characteristics of the FOD operator, a suitable fractional order v is selected.
(4) The amplification factor k of the fractional-order differential operator is calculated based on the fractional order v and the data function G(x).
(5) The energy loss of the signal transmitted under the existing conditions is calculated by formulating the transmission method based on the characteristics of the signal G(t).
(6) The amplification factor Kg required for the effective transmission of signal G(t) is calculated by the detection system based on the transmission distance L of the system.
(7) The required amplification factors Kg and K are combined, and the number of differential processing cycles m required by the system is calculated to realize the remote transmission target of information data set by the underground space detection system.

# 5 THE QUALITY OF INFORMATION DATA IMPROVING METHOD FOR UNDER-GROUND SPACE INFORMATION BASED ON FOD ALGORITHM

## 5.1 Fundamental Principle

The quality of information data improving is essentially a dual function achieved through the processing of the detection data. As mentioned above, the transmission of data over a desired distance can be achieved by setting the number of iterations *j* of the differential operator for data processing based on the amplification factor *K* of the FOD operator and the established amplification factor $K_G$ of the detection system. The detection accuracy of the data can be achieved by comparing the established accuracy $S_G$ of the system and the accuracy *S* of the FOD operator and gradually adjusting the step value *h* of the FOD operator so that $S_G$ can be approximately achieved. Therefore, obtaining the step value *h* and the amplification factor *K* required by the FOD operator in data processing under the established conditions is a prerequisite for achieving the research objectives. According to the results of the above research, it is known that the accuracy of the data detection at the same fractional order depends on the step value *h* of the FOD operator, independent of the amplification factor *K*. The calculation of the amplification factor *k* of the FOD operator must be based on the fact that the step value *h* is known. Therefore, we can improve the quality of underground spatial information data by sequentially calculating the step size value *h* required to meet the established goal, and then applying it to the method of calculating the amplification coefficient *K* of the FOD operator

## 5.2 Mathematical Model

The analysis system is intended to transmit an energy signal *H(t)* to a location at a distance of *L*, and the accuracy of the acquired data must not be lower than $S_G$. To ensure the reliability of the data detection system, the detection signal needs to be amplified, and the enhancement factor must not be lower than $K_G$. The data $H(x_i)$ of the signal *H(t)* at different time points $t_i$ are collected. The main impact parameter of the obtained data is *x*, and its value interval [*a, b*]. A function *H(x)* relating the detection data value $H(x_i)$ and $x_i$ is fitted. According to the signal characteristics of the signal *H(t)* and the spectral characteristics of the FOD operator, the proposed method of applying the differential operator at fractional order *v* for the long-range, high-precision detection of data is used. With reference (1), the FOD processing equation of data under the known conditions is

$$_aD_x^v H(x) = \lim_{h \to 0} H_h^v(x) = h^{-v} \sum_{i=0}^{n} (-1)^i \binom{v}{i} H(x-ih). \tag{21}$$

In (21), $\binom{v}{i}$ is a binomial coefficient, defined as

$$\binom{v}{i} = \frac{v(v-1)(v-2)....(v-i+1)}{i!}. \tag{22}$$

According to the above theory for the long-distance and high-precision transmission of data, the standard deviation of the data is calculated after the FOD operator is applied as follows:

$$S = \left( \sum_{i=0}^{n} (H_{a+ih}^v - \bar{H}^v) / K \right)^{0.5}. \tag{23}$$

If the detection accuracy target set by the system is to be achieved, the following conditions must be met:

$$\begin{cases} S \leq S_G \\ S_G - S \approx 0 \end{cases}. \tag{24}$$

As mentioned above, if the step value of the FOD operator is *h*, the system can be calculated by a stepwise approximation method to meet the system accuracy target $S_G$. Then, *h* is introduced to (23) to obtain the amplification factor of the detection data after processing by the *v*-order differential operator as follows:

$$K = \bar{H}^v(x) / \bar{H}(x) = h^{-v} \sum_{i=0}^{n} \sum_{i=0}^{n} (-1)^i \binom{v}{i} H(x-ih) / \sum_{i=0}^{n} H(a+ih). \tag{25}$$

By combining the required data enhancement factor $K_G$ of the detection system and referring to (19), the number of fractional-order differentiations *m* required to detect the data is

$$m = K_G / K = K_G \bar{H}(x) / \bar{H}^v(x) = K_G h^v \sum_{i=0}^{n} H(a+ih) / \sum_{i=0}^{n} \sum_{i=0}^{n} (-1)^i \binom{v}{i} H(x-ih). \tag{26}$$

The number of fractional order differentiation processes for the detection data is *m*.

## 5.3 Implementation Steps

Suppose the underground space information detection system detects an energy signal *H(t)*, where $t \in (t_1, t_2)$. Referring to (7), we can get the corresponding equation *G(x)*, $x \in (a, b)$. The signal characteristics of the data to be detected and its usage requirements are analyzed. The data accuracy $S_G$ and the amplification factor $K_G$ required to achieve the data detection target are set. According to the above fundamental principle, The remote and high-precision detection methods of underground space information data are as follows:

(1) A data processing scheme is developed for underground space information data based on fractional calculus theory, according to the set data accuracy SG and signal amplification coefficient KG.
(2) The collected data H(xi) of underground space information is applied to analyze the important impact parameter x, and the functional relationship between the data and the impact parameter H(x) is obtained by fitting.

(3) The appropriate fractional order v for the processing of data is selected based on the characteristics of the detection data signal and the amplitude and frequency characteristics of the FOD operator.
  (4) The step value h required to achieve the required accuracy SG of the system is calculated based on the mathematical model shown in (12) and (14).
  (5) The amplification factor k of the FOD operator is calculated based on the data, and the number of differential processes $m$ required to achieve the data amplification factor $K_G$ is determined.

Through the integration of high-precision detection method and long-distance transmission method of underground space information data, the quality of information data improving function of a signal in underground space is realized. It not only successfully solves the problem that the two functions of the current method cannot be considered, but also realizes the detection technology of information data under the established detection target.

## 6 APPLICATION EXAMPLES

### 6.1 Experimental Environment

Huainan City is an important coal production base in China, and the prominent feature of its coal mines is that there are more "high-gas" mines. The output of the high-gas mining faces accounts for more than 70% of the total output of the region. Therefore, the control of the gas concentration at the working face has been a key issue and difficulty for ensuring safe production in the mine. To realize the effective control of the gas concentration at the working surface, we first need to accurately measure the gas concentration at the working surface and transmit the detection data in real time. However, the transmission distance of the detection data and the detection accuracy of the data are limited by the properties of the detection instruments, working environment, signal interference, and other factors. Therefore, achieving accurate measurement and long-distance transmission of the gas concentration at the working surface has always been a technical problem in the safety management of high-gas coal mines.

To test the application of fractional-order calculus theory in data detection, the 151302 working surface of a mine in Huainan city was considered, which is located about 500 m from the surface and about 2 km from the main shaft. The test required real-time accurate detection of the gas concentration at the working site under the existing conditions. The experiment required the real-time accurate gas concentration detection function for the operation site under the existing conditions. Due to the long distances and dynamic changes in the operation site at the coal mine working face, the experiment required a wireless network to realize real-time long-distance transmission of the collected data. To ensure the normal production of the working surface and the safety of the experiment site, after considering various factors, we designed a structure diagram of the gas concentration detection system, which is shown in Figure 4, and developed the corresponding data detection scheme. The data sensed by the gas concentration sensor was first transmitted to the gas concentration detection system near the working face by Wi-Fi wireless network transmission, then the data processed by fractional-order differentiation was transmitted to the underground data center by an optical fiber transmission method. Finally, the data collected by the underground data center was transmitted to the surface data management center by Controller Area Network (CAN).

Due to the complex spatial environment of the underground structure and the twists and turns of the fiber arrangement path, the traditional calculation of fiber energy attenuation coefficient is not suitable for this case. In order to obtain accurate experimental results, this paper adopts the experimental method to calculate the energy attenuation coefficient during the signal transmission under the existing working conditions. In the experiment, we choose the information transmission distance of 100 meters, and found its energy attenuation coefficient is 0.63 after testing. Therefore, the data management center needed to amplify the signal strength by more than 6.25 times to achieve effective sensing and high-precision collection of data on the operating surface 2 km away. To ensure the detection quality of the collected data, the standard deviation of the data should not be greater than 0.005.

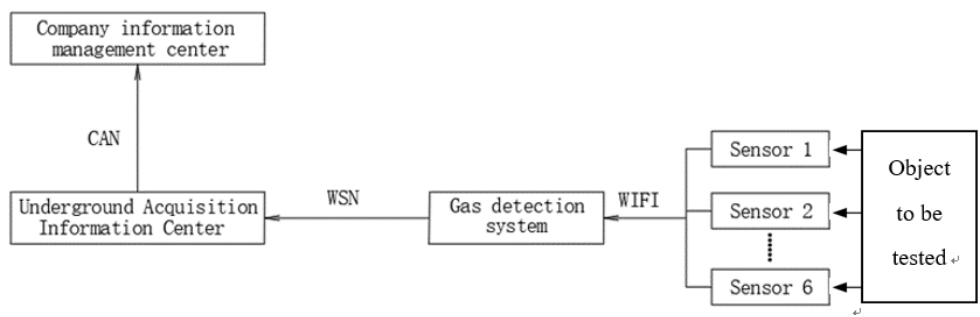

Fig. 4. Structure diagram of gas concentration detection system

To improve the accuracy of the data, six working 3-2.3V LEL combustible gas sensors were used to measure the gas concentration at the working site of the working face. The six sensors worked simultaneously and measured once every 5 s, and the data collection center collected five gas concentration measurements from six sensors. Therefore, we can extract 30 data at the same location in real time as experimental samples.as shown in Table 1.

Table 1. Experimental Data Table (%)

| Sensor No. | 1# | 2# | 3# | 4# | 5# | 6# |
|---|---|---|---|---|---|---|
| 1st measured value | 10.28 | 11.02 | 11.06 | 10.48 | 11.16 | 10.62 |
| 2nd measured value | 10.85 | 10.35 | 10.66 | 10.25 | 10.45 | 10.55 |
| 3rd measured value | 10.55 | 10.28 | 11.25 | 10.86 | 11.35 | 10.36 |
| 4th measured value | 11.08 | 11.31 | 11.35 | 11.08 | 10.58 | 11.28 |
| 5th measured value | 10.75 | 10.56 | 10.58 | 10.45 | 10.25 | 10.55 |

## 6.2 Pre-processing of Data

1) ANALYSIS OF DATA

In the experiment, the average value of five measurements of each sensor, $\bar{E_i}$, was considered to be the true value, and the average value of all the measurements, $\bar{E}$, was considered to be the true value, $E$. The summary of the gas concentration detection data of the working face is shown in Table 2. The detection data showed an irregular distribution around the measured true value. Therefore, influenced by multiple factors, such as the properties of the measuring instrument FOD and the working environment, the initial detection value of the gas concentration at the coal mine working surface had a large detection error, which would create difficulties in providing accurate detection values for the safety management of the coal mine and significantly affect the decision making of the safety management system.

Table 2. Summary of Experimental Data (%)

| Sensor No. | 1# | 2# | 3# | 4# | 5# | 6# |
|---|---|---|---|---|---|---|
| Average value $\bar{E_i}$ | 10.702 | 10.704 | 10.980 | 10.624 | 10.758 | 10.672 |
| Standard deviation $S$ | 0.27 | 0.40 | 0.31 | 0.30 | 0.42 | 0.32 |
| True value $E$ | 10.740 | | | | | |
| System standard deviation $S_x$ | 0.1146 | | | | | |

2) IMPACT PARAMETER OF DATA

Due to the harsh environment and complex working conditions of the underground structure space, the detection data in the transmission process will inevitably be affected by a variety of factors, which bringing a large data detection error. All sensors are placed in the same place, and they operated at the same time and with the same frequency to detect the gas concentration at the location. Thus, the effects of the environment and energy losses in the transmission of data were basically the same, this means that the property of measuring instruments directly affects the measured value of information data. Since the standard deviation can reflect the property of measuring instrument effectively, it was used as the impact parameter of the measured values of the sensor to analyze the correlation between the property of the testing equipment and the measured value.

3) EQUATION RELATING CONCENTRATION DATA AND IMPACT FACTOR

Since least square algorithm has good data fitting accuracy. Therefore, this paper applies it to the fitting of the distribution function of underground spatial information data. According to the parameter values $E_i$ and $x_i$ shown in Table 2, the mathematical expression of the function $E(x)$ was assumed to be

$$E(x) = a_0 + a_1 x + a_2 x^2 + \cdots + a_n x^n . \qquad (27)$$

According to formula (27), we can know that the function formula $E(x)$ depends on the parameters $n$ and $a_i$ ($i = 0, 1,... n$) in the formula. For the characteristics of least square algorithm, the parameter $n$ of function $E(x)$ should not exceed the number of samples, so $n \in [0, 6]$, and there is an optimum was found to minimize the error between the function $E(x)$ and the measured true value $E$ when $n \in [0, 6]$. The Polyfit function in the MATLAB software was applied, and the total errors of the fitted values of (27) at different orders are shown in Table3. From the simulation results, we can see that when

$n = 3$, the error between the fitting value and the measured true value is the smallest, and the value is 0.012. At this time, the equation of the $E(x)$ function was

$$E(x) = 0.035x - 10.727.  \quad (28)$$

Table 3. Summary of Experimental Data (%)

| Order $n$. | 1 | 2 | 3 | 4 | 5 |
|---|---|---|---|---|---|
| Total error | 0.012 | 0.096 | 0.174 | 0.178 | 0.520 |

### 6.3 High-precision Detection of Data Based on Stepwise Approximation Method

According to the working principle of the stepwise approximation method and the characteristics of the data in this case, the implementation process of the method was as follows.
  (1) The corresponding mathematical treatment under the FOD was modeled based on the available experimental conditions and data.
  (2) The appropriate differential order v was selected based on the detection data characteristics and the FOD operator properties.
  (3) The established mathematical model was applied to set the initial step value h, and the initial standard deviation Ss was calculated under known conditions.
  (4) The step value h was continuously adjusted by comparing with the set system threshold SG until the accuracy met the set system accuracy threshold $S_G$.

1) MATHMATICL MODEL OF DATA ACCURACY CALCUTION BASED ON FOD OPERATOR

According to (4) and Figure 1, when the differential order $v \in [0, 1]$, the signal strength increases with the increase of the fractional order $v$ at the high frequency stage, but with the increase of the frequency, the difference between the enhancement values of the FOD operator on the signal at different orders shows a downward trend. Therefore, from the perspective of space saving and convenient calculation, this case discusses the application effect of FOD operator in underground space information data fusion processing when the median value between fractional orders [0, 1] is $v = 0.5$. According to (7), we can get the formula for the standard deviation of the data in this case.

$$S = \left( \frac{\sum_{i=0}^{[b-a]/h} \left( E^{0.5}(x+ih) - \bar{E}^{0.5}(x) \right)^2}{nK} \right)^{0.5}, \quad (29)$$

where

$$E^{0.5}(x+ih) = \lim_{h_4 \to 0} E_h^{0.5}(x+ih) = h^{-0.5} \sum_{i=0}^{[b-a]/h} (-1)^i \binom{0.5}{i} E(x-ih). \quad (30)$$

Combined with (28), the data shown in the table were processed by the 0.5-order differential operator, and the result was

$$E^{0.5}(x_i) = h^{-0.5} \left( \begin{array}{c} E(x_i) - 0.5E(x_i - h) - \dfrac{0.5(1-0.5)}{2!} E(x_i - 2h) - \text{L} \\ - \dfrac{\Gamma(1-0.5)}{n!\Gamma(n+1-0.5)} E(x_i - nh) \end{array} \right), \quad (31)$$

where $b = 0.26$; $a = 0.07$; $n = [b-a]/h = [0.26 - 0.07]/h = 0.19$.

2) CALCULATION OF STEP SIZE $H$ BASED ON STEP-STEP APPROACH METHOD

To improve the efficiency of the system, a step-by-step method was used to adjust the step value $h$ to achieve continuous improvement in the accuracy of the system data and select the best step value in terms of data accuracy and efficiency. Apply MATLAB mathematical simulation software to calculate the detection accuracy of the detection system was calculated by initially selecting $h = 0.01$. The data processing model for the detection data with step $h = 0.01$ and order $v = 0.5$ was obtained according to (29) and (31):

$$E_{0.01}^{0.5}(x) = 10 \sum_{i=0}^{19} (-1)^i \binom{0.5}{i} E_{0.01}(x - 0.19). \quad (32)$$

According to the data shown in Table 2 and (32), the experimental data shown in Table 4 were obtained for fractional order $v = 0.5$ and the step size $h = 0.01$.

Table 4. Experimental data summary table when $h = 0.01$ (%)

| Sensor No. | 1# | 2# | 3# | 4# | 5# | 6# |
|---|---|---|---|---|---|---|

| Average value $\bar{E}_i$ | 10.702 | 10.904 | 10.98 | 10.624 | 10.758 | 10.672 |
|---|---|---|---|---|---|---|
| Standard deviation $S_i$ | 0.27 | 0.40 | 0.31 | 0.30 | 0.42 | 0.32 |
| Measured true value $E$ | 10.7400 | | | | | |
| Value after fusion $E_{0.01}^i$ | 15.524 | 15.532 | 15.526 | 15.526 | 15.533 | 15.527 |
| Post-fusion mean $\bar{E}_{0.01}$ | 15.5279 | | | | | |
| Amplification factor $K_{0.01}$ | 1.446 | | | | | |

To facilitate the comparison of the standard deviation between the data before and after data processing, the processing results shown in Table 4 were divided by the amplification factor $K_{0.01}$, and the results are shown in Table 5.

Table 5. Experimental data summary table when $h$ = 0.01 (%)

| Sensor No. | 1# | 2# | 3# | 4# | 5# | 6# |
|---|---|---|---|---|---|---|
| Average value $\bar{E}_i$ | 10.702 | 10.904 | 10.98 | 10.624 | 10.758 | 10.672 |
| Standard deviation $S_i$ | 0.27 | 0.40 | 0.31 | 0.30 | 0.42 | 0.32 |
| Fusion final value $\bar{E}_{0.01}^i$ | 10.736 | 10.741 | 10.737 | 10.737 | 10.742 | 10.738 |
| Pre-fusion standard deviation $S$ | 0.1146 | | | | | |
| Post-fusion standard deviation $S_p$ | 0.0027 | | | | | |

When $h$ = 0.01, the accuracy of the fused data was significantly lower than threshold value of 0.005. Based on the error between the two, the step value was adjusted to $h$ = 0.005. From (31), we determined that the data processing model at this time is

$$E_{0.005}^{0.5}(x) = 14.144 \sum_{i=0}^{38}(-1)^i \binom{0.5}{i} E_{0.005}(x - 0.19) \cdot \quad (33)$$

Based on the data shown in Table 2 and (33), the processing results shown in Table 6 can be obtained.

To compare the standard deviation of data before and after data processing, the data processing results shown in Table 6 were divided by the amplification factor $K_{0.005}$ to obtain the data value shown in Table 7.

Table 6. Experimental data summary table when $h$ = 0.005 (%)

| Sensor No. | 1# | 2# | 3# | 4# | 5# | 6# |
|---|---|---|---|---|---|---|
| Average value $\bar{E}_i$ | 10.702 | 10.904 | 10.98 | 10.624 | 10.758 | 10.672 |
| Standard deviation $S_i$ | 0.27 | 0.40 | 0.31 | 0.30 | 0.42 | 0.32 |
| Post-fusion mean $E$ | 10.7400 | | | | | |
| Value after fusion $E_{0.005}^i$ | 15.865 | 15.873 | 15.867 | 15.866 | 15.874 | 15.868 |
| Mean after fusion $\bar{E}_{0.005}$ | 15.8387 | | | | | |
| Amplification factor $K_{0.005}$ | 1.4552 | | | | | |

Table 7. Experimental Data Summary Table When $h$ = 0.01 (%)

| Sensor No. | 1# | 2# | 3# | 4# | 5# | 6# |
|---|---|---|---|---|---|---|
| Average value $\bar{E}_i$ | 10.702 | 10.904 | 10.98 | 10.624 | 10.758 | 10.672 |
| Standard deviation $S_i$ | 0.27 | 0.40 | 0.31 | 0.30 | 0.42 | 0.32 |
| Fusion final value $\bar{E}_{0.005}^i$ | 10.738 | 10.743 | 10.739 | 10.739 | 10.744 | 10.740 |
| Pre-fusion standard deviation $S$ | 0.1146 | | | | | |
| Post-fusion standard deviation $S_p$ | 0.0021 | | | | | |

When $h = 0.005$, the accuracy between the fused data was slightly lower than the threshold value, which was selected to improve the efficiency. Next, $h = 0.003$ was selected, and we obtained the following processing model of the data at this time:

$$E_{0.003}^{0.5}(x) = 18.257 \sum_{i=0}^{63} (-1)^i \binom{0.5}{i} E_{0.003}(x-0.19). \tag{34}$$

To compare the standard deviation of the data before and after data processing, the data processing results shown in Table 4 were divided by the amplification factor $K_{0.003}$ to obtain the values shown in Table 8.

The data obtained by this process met the system accuracy requirement after being processed by the differential operator with fractional order $v = 0.5$ and step value $h = 0.003$.

Table 8. Experimental Data Table When $h = 0.003$ (%)

| Sensor No. | 1# | 2# | 3# | 4# | 5# | 6# |
|---|---|---|---|---|---|---|
| Average value $\bar{E}_i$ | 10.702 | 10.904 | 10.98 | 10.624 | 10.758 | 10.672 |
| Standard deviation $S_i$ | 0.27 | 0.40 | 0.31 | 0.30 | 0.42 | 0.32 |
| Measured true value $E$ | 10.7400 | | | | | |
| Value after fusion $E_{0.003}^i$ | 16.063 | 16.068 | 16.064 | 16.064 | 16.071 | 16.065 |
| Post-fusion mean $\bar{E}_{0.003}$ | 16.0659 | | | | | |
| Amplification factor $K_{0.003}$ | 1.49589 | | | | | |

## 6.4 Long-distance data transmission based on cyclic iteration method

### 1) Value of parameter $h$ and amplification factor $k$

Since the strength of the signal depended on its step value $h$ and the number of FOD processing iterations, the original data was next processed for the first iteration through the differential results with $h = 0.003$, $v = 0.5$, and $x \in [0.27, 0.42]$. Before and after the iteration, the applied parameters were the same in the mathematical model. Various data points before the iterations are shown in Table 4.

By combining the FOD operator shown in (35) with the data processing model, the data shown in Table 9 were processed by the first FOD operator iteration to obtain the values shown in Table 10. The amplification coefficient and data accuracy before the iterations were basically equal between the data after iterative processing by the FOD operator, thus verifying the correctness of the above assertion about the relationship between the accuracy and amplification coefficient in the data processing and the FOD parameters. When $v = 0.5$ and $h = 0.003$, the amplification factor $k \approx 1.47$ for the data before and after the iteration of the FOD operator.

Table 9. Experimental Data Bable When $h = 0.003$ (%)

| Sensor No. | 1# | 2# | 3# | 4# | 5# | 6# |
|---|---|---|---|---|---|---|
| Average value $\bar{E}_i$ | 10.702 | 10.904 | 10.940 | 10.624 | 16.071 | 16.065 |
| Standard deviation $S_i$ | 0.27 | 0.40 | 0.31 | 0.30 | 0.42 | 0.32 |
| Fusion final value $\bar{E}_{0.003}^i$ | 10.738 | 10.741 | 10.739 | 10.739 | 10.744 | 10.739 |
| Pre-fusion standard deviation $S$ | 0.1146 | | | | | |
| Post-fusion standard deviation $S$ | 0.0019 | | | | | |

Table 10. Summary of Experimental Data (%)

| Sensor No. | 1# | 2# | 3# | 4# | 5# | 6# |
|---|---|---|---|---|---|---|
| Pre-iterative data | 10.738 | 10.741 | 10.739 | 10.739 | 10.744 | 10.739 |
| Pre-iteration mean $E$ | 10.74 | | | | | |
| Standard deviation $S$ | 0.27 | 0.40 | 0.31 | 0.30 | 0.42 | 0.32 |
| Data after iteration $\bar{E}_d^i$ | 16.061 | 16.068 | 16.063 | 16.063 | 16.069 | 16.064 |
| Standard deviation after iteration $S$ | 0.00188 | | | | | |

| Amplification factor $k$ | 1.46575 |
|---|---|

**2) Values of amplification factor $K$ and number of iterations $m$**

From the above discussion on fractional-order calculus theory, the data processing accuracy, and the amplification coefficient, combined with the results for this example, the amplification coefficient $k$ of the FOD operator on the data was obtained based on the relationship between the data accuracy $S$ and the number of iterations. The relationships were as follows:

$$\begin{cases} k = 1.47 \\ S = 0.0019 \end{cases}, \quad (35)$$

If $q$ is the number of FOD processing iterations of the data. For the detection system given and (19), to meet the amplification factor set for the detection system, it was necessary to satisfies the following equations:

$$\begin{cases} K = 1.47^q \\ 6.25 \leq 1.47^q < 1.47^{q+1} \end{cases} \quad (36)$$

By solving for $q$ in (36), we can see that under the above detection conditions, when the FOD operator step $h = 0.003$, the number of FOD processing cycles $q = 5$. The amplification coefficient of the detection signal $K = 6.86$, which meets the set system requirements.

### 6.5 Experimental Results and Their Analysis

According to the above research results, it can be concluded that: when the fractional order $v = 0.5$ and the step size $h = 0.003$, the detection accuracy of the information data in Table 1 is 0.0019 after the experimental data is processed by fractional order differentiation as shown in Table 9. the amplification coefficient of the information data is 6.86 after 5 times processing by fractional order differentiation operator. The final information data processing results are shown in Table 11, from which we could see that the detected information data can meet the detection target set by the system completely.

Table 11. Summary of Final Processing Results of Experimental Data (%)

| Sensor No. | 1# | 2# | 3# | 4# | 5# | 6# |
|---|---|---|---|---|---|---|
| Average value $\bar{E}_i$ | 10.702 | 10.904 | 10.94 | 10.624 | 16.071 | 16.065 |
| Standard deviation $\sigma_i^2$ | 0.27 | 0.40 | 0.31 | 0.30 | 0.42 | 0.32 |
| Fusion final value $\bar{E}_{0.003}^i$ | 10.736 | 10.741 | 10.737 | 10.737 | 10.742 | 10.738 |
| Final detection value $E_z$ | 73.649 | 73.683 | 73.656 | 73.656 | 73.690 | 73.663 |
| Standard deviation before treatment $S$ | 0.1146 | Standard deviation after treatment $S_4$ | | | | 0.0019 |
| Final amplification factor $K$ | 6.86 | | | | | |

Since the experimental site in the paper was located in a structural space more than 500m underground, it was a harsh experimental environment. The experiments were conducted by simplifying the calculation steps and selecting typical parameter values to analyze and process the experimental data. In practical applications, the corresponding fractional order $v$ can be selected based on the data characteristics, and the detection system can achieve the detection distances and data accuracies arbitrarily by adjusting the step value $h$ and processing number $n$ of the FOD operator.

## 7 CONCLUSION

Fractional-order calculus theory was applied to the detection of underground space information data. by analyzing the characteristics of the generate information in complex environments and the signal amplitude and frequency characteristics after fractional differentiation algorithm processing to meet the needs of improving the quality of information data. The fractional calculus theory is successfully applied to the field of information data detection. The application experiments prove that the method described in this paper has the following characteristics:

(1) According to the characteristics of FOD operator, the real-time detection function of various information data in underground space can be realized.
(2) By adjusting the step value h, the high-precision detection technology of information data in the complex environment of underground space is realized.
(3) By adjusting the number n of FOD processing, the remote detection technology of underground space in-formation data is realized.

(4) Through the integration of high-precision detection method and long-distance transmission method, the high-precision remote transmission function of a signal in underground space is realized.
(5) Using the integration of high-precision detection method and long-distance transmission method, the underground space information detection system can realize the information data detection technology of setting property goals.

Therefore, the algorithm used in this paper can successfully solve the problems existing in the current underground space information detection methods.

**Acknowledgements:** This project was supported by the National Natural Science Foundation, China [grant numbers 51878005 and 51778004] and the Anhui Provincial department of Education [grant number KJ2020A0488].